\documentclass{amsart}

\usepackage[utf8]{inputenc}
\usepackage{mathrsfs, mathtools}
\usepackage{amsfonts, amssymb, amsmath, amsthm}
\usepackage[pagebackref]{hyperref}
\hypersetup{final, colorlinks=true, citecolor=blue}
\usepackage{url}
\usepackage{enumitem}
\usepackage{caption}
\usepackage{stmaryrd}
\usepackage{xcolor}
\usepackage{orcidlink}

\newcommand{\Q}{\mathbb{Q}}

\newcommand{\Z}{\mathbb{Z}}

\newcommand{\pb}{\overline{p}}
\newcommand{\eml}{\epsilon_{m,\ell}}
\newcommand{\Zq}{\Z\llbracket q\rrbracket}
\newcommand{\Um}[1][m]{\vert U(#1)}
\newcommand{\Vm}[1][m]{\vert V(#1)}
\newcommand{\Tls}{\vert T(\ell^2)}

\newcommand{\eigen}{\lambda_{m,\ell}}
\newcommand{\calV}{\mathcal{V}}

\DeclareMathOperator{\SL}{SL}

\DeclareMathOperator{\val}{val}
\DeclareMathOperator{\sgn}{sgn}

\DeclareMathOperator{\ord}{ord}

\def\ll#1{{\left\langle{#1}\right\rangle}}
\def\smat#1{\left(\begin{smallmatrix} #1 \end{smallmatrix} \right)}
\def\pmat#1{\left(\begin{matrix} #1 \end{matrix} \right)}

\newcommand{\set}[2]{\left\{{#1} \, : \, {#2}\right\}}
\newcommand{\kro}[2]{\genfrac(){}{}{#1}{#2}}

\newcommand{\modf}[2]{\mathcal{M}_{#1}({#2})}
\newcommand{\cusf}[2]{\mathcal{S}_{#1}({#2})}
\newcommand{\eisf}[2]{\mathcal{E}_{#1}({#2})}

\newcommand{\modfc}[3]{\mathcal{M}_{#1}({#2,#3})}
\newcommand{\cusfc}[3]{\mathcal{S}_{#1}({#2,#3})}
\newcommand{\eisfc}[3]{\mathcal{E}_{#1}({#2,#3})}

\numberwithin{equation}{section}

\theoremstyle{plain}
\newtheorem{prop}[equation]{Proposition}
\newtheorem{thm}[equation]{Theorem}
\newtheorem{coro}[equation]{Corollary}
\newtheorem{lemma}[equation]{Lemma}

\theoremstyle{remark}

\newtheorem{rmk}[equation]{Remark}

\begin{document}

\title
{Explicit families of congruences for the overpartition function}

\author
[
	N. C. Ryan,
	N. Sirolli,
	J.C. Villegas-Morales and
	Q.-Y. Zheng
]
{
	Nathan C. Ryan\orcidlink{0000-0003-4947-586X},
	Nicolás Sirolli\orcidlink{0000-0002-0603-4784},
	Jean Carlos Villegas-Morales\orcidlink{0009-0004-8933-0627} and
	Qi-Yang Zheng
}

\address{Department of Mathematics, Bucknell University, Lewisburg, Pennsylvania 17837} 
\email{nathan.ryan@bucknell.edu}

\address{Departamento de Matemática - FCEyN - UBA and IMAS - CONICET, Pabellón I, Ciudad Universitaria, Ciudad Autónoma de Buenos Aires (1428), Argentina}
\email{nsirolli@dm.uba.ar}

\address{Escuela de Matem\'atica, Universidad de Costa Rica, San Jos\'e 11501, Costa Rica}
\email{jean.villegas@ucr.ac.cr}

\address{Department of Mathematics, Sun Yat-sen University (Zhuhai Campus), China}
\email{zhengqy29@mail2.sysu.edu.cn}

\begin{abstract}

	In this article we exhibit new explicit families of congruences for the
	overpartition function, making effective the existence results given
	previously by Treneer.
	We give infinite families of congruences modulo $m$ for $m = 3, 5, 7, 11$,
	and finite families for $m = 13, 17, 19$.
	
\end{abstract}

\maketitle

\section{Introduction}

Let $p(n)$ be the number of partitions of a positive integer $n$; that is, the
number of ways $n$ can be written as a sum of non-increasing positive integers.
Ramanujan \cite{ramanujan} proved congruences of the form:
\begin{align*}
p(5n+4) & \equiv 0 \pmod{5},\\
p(7n+5) & \equiv 0 \pmod{7},\\
p(11n+6) &\equiv 0 \pmod{11},
\end{align*}
for every $n$.
For decades it was difficult to find more congruences like these; nevertheless,
Ono proved in \cite{ono} that for each prime $m \geq 5$ there exists an infinite
family of congruences for the partition function modulo $m$: more precisely, he proved
that a positive proportion of the primes $\ell$ are such that
\begin{align*}
	p\left(\frac{m\ell^3 n + 1}{24}\right) \equiv 0  \pmod m.
\end{align*}
for every $n$ prime to $\ell$.

\medskip 

The number of overpartitions $\overline{p}(n)$ of a positive integer $n$ is
defined to be the number of ways in which $n$ can be written as a non-increasing
sum of positive integers in which the first occurrence of a number may be
overlined (see \cite{overpartitions}).

The numbers of both partitions and overpartitions can be described in terms of
eta-quotients; in particular, they are known to be coefficients of weakly
holomorphic modular forms of half-integral weight, with integral coefficients.
Treneer showed in \cite{treneer2} that Ono's existence results were
valid, more generally, for the coefficients of such modular forms.
In the particular case of the overpartition function 
we improve her result obtaining the following theorem.

\begin{thm}
	\label{thm:treneer}
	Let $m$ be an odd prime.
	Then a positive proportion of the primes $\ell \equiv -1 \pmod{16m}$
	have the property that
	\begin{align*}
		\overline p\left(m\ell^3 n\right) \equiv 0  \pmod m.
	\end{align*}
	for every $n$ prime to $m\ell$.
\end{thm}

The main goal of this article is to exhibit explicit instances of these (families
of) congruences, as well as for certain variations similar to those considered
by Ono for the partition function.

\medskip

Weaver devised a strategy in \cite{weaver} for making Ono's results explicit:
she exhibited 76,065 new families of congruences for the partition
function by finding congruences between its generating function and appropriate holomorphic modular forms, and then
verifying a finite number of congruences for the partition function.
Her computations were extended by Johansson \cite{johansson}, who used efficient
algorithms for computing the partition function to find
more than $2.2 \cdot 10^{10}$ such families of congruences.

Using Weaver's techniques along with the theory of Eisenstein series of
half-integral weight from \cite{wang-pei}, we were able to find \emph{infinitely
many} families of congruences for the overpartition function.
Our first main results are the following two theorems.

For an odd prime $m$, throughout the article we denote
\begin{equation*}
	\label{eqn:km}
	k_m =
	\begin{cases}
		m + 2, & m = 3 \\
		m - 2, & m > 3.
	\end{cases}
\end{equation*}

\begin{thm}
	\label{thm:inf0}
	Let $m \in \{3,5,7,11\}$, and let $\ell$ be an odd prime such that
	$\ell^{k_m-2} \equiv -1 \pmod m$.
	Then
	\begin{equation*}
		\pb\left(m \ell^3 n\right) \equiv 0 \pmod m
	\end{equation*}
	for every $n$ prime to $\ell$.

\end{thm}

We remark that for $m=3$ and $m=5$ the result was proved, respectively, in
\cite[Coro. 1.5]{lovejoy2011quadratic} and \cite[Prop. 1.4]{treneer1}.
We include those cases in our results to highlight our unified approach.

\begin{thm}
	\label{thm:infeps}
	Let $m \in \{3,5,7,11\}$, and let $\ell$ be an odd prime such that
	$\ell^{k_m-2} \equiv -1 + \eml \, \ell^{\tfrac{k_m-3}2} \pmod m$,
	with $\eml \in \{\pm1\}$.
	Then
	\begin{equation*}
		\pb\left(m \ell^2 n\right) \equiv 0 \pmod m
	\end{equation*}
	for every $n$ prime to $\ell$ such that
	\begin{equation*}
		\kro{(-1)^{\tfrac{k_m-1}2}n}{\ell} = \eml.
	\end{equation*}
\end{thm}

For primes $m \geq 13$ the appearance of cuspidal forms in level $16$ and weight
$k_m/2$ makes it more difficult to find infinitely many families of
congruences.
Using the results from \cite{barquero-etal} for efficiently computing the
overpartition function, we obtain the following families of congruences.

\begin{thm}
	\label{thm:fin0}
	Let $m,\ell$ be primes as in Table~\ref{tab:ml}.
	Then
	\begin{equation*}
		\pb\left(m \ell^3 n\right) \equiv 0 \pmod m
	\end{equation*}
	for every $n$ prime to $\ell$.
\end{thm}

\begin{table}[ht]
\begin{tabular}{c|l}
	$m$ & $\ell$
    \\\hline
	$13$ & $1811, 1871, 1949, 2207, 3301, 4001, 4079, 4289, 4931$ \\ 
	$17$ & $2039, 2719, 3331, 4079$ \\ 
	$19$ & $151, 1091, 2659, 3989$ \\
\end{tabular}
\caption{Congruences for primes $m \geq 13$.  See Theorem~\ref{thm:fin0}.}
\label{tab:ml}
\end{table}

\begin{thm}
	\label{thm:fineps}
	Let $m,\ell$ be primes, and let $\eml \in \{\pm1\}$
	be as in Table~\ref{tab:ml_eps}.
	Then
	\begin{equation*}
		\pb\left(m \ell^2 n\right) \equiv 0 \pmod m
	\end{equation*}
	for every $n$ prime to $\ell$ such that
	\begin{equation*}
		\kro{(-1)^{\tfrac{k_m-1}2}n}{\ell} = \eml.
	\end{equation*}
\end{thm}

\begin{table}[ht]
\begin{tabular}{c|l}
	$m$ & $(\ell,\eml)$
    \\\hline
	$13$ & $(431, 1), (2459, 1), (4513, 1), (4799, 1)$ \\
	$17$ & $(167, 1), (541, 1), (911, -1), (1013, -1), (1153, 1)$,\\
		 & $(1867, 1), (1931, -1),(2543,-1), (2683, 1), (2887, 1)$,\\
		 & $(3019, -1), (3023, 1), (3329, 1), (4243, -1), (4651, -1)$ \\
	$19$ & $(2207,-1)$ \\
\end{tabular}
\caption{Congruences for primes $m \geq 13$.  See Theorem~\ref{thm:fineps}.}
\label{tab:ml_eps}
\end{table}

We point out that using different techniques, in \cite{rsst,barquero-etal} the
authors found (finite) families of congruences for the overpartition function
modulo $m$ for $m = 3,5,7$;  
see also \cite{Chen} for $m = 5$, and \cite{Xia} for powers of $m=3$.
As far as we know, the results in this article give the first known congruences
for $m > 7$.

\bigskip

The rest of the article is organized as follows.
In the next section we give the necessary notation and preliminaries regarding
half-integral weight modular forms and eta-quotients.
In Section~\ref{sect:eisen} we state the results we need on Eisenstein series of
half-integral weight and level 16.
We conclude the article with the proofs of our main results in
Section~\ref{sect:proofs}.

\section{Preliminaries} \label{sect:prelims}

\subsection*{Half-integral weight modular forms}

We refer the reader to \cite[Sect. 5]{wang-pei} for details on this subsection.

\medskip

Given a non zero integer $m$ we denote by $\chi_m$ the primitive Dirichlet
character such that $\chi_m(a) = \kro ma$ for every $a$ such that $(a,4m) = 1$.

Given an odd integer $k\geq 3$, we denote $\lambda = \tfrac{k-1}2$.
Furthermore, given a positive integer $m$ we denote 
$\omega_n = \chi_m$, with $m = (-1)^\lambda n$.

Given $k$ as above, a positive integer $N$ divisible by $4$ and a
character $\chi$ modulo $N$,
we denote by 
$\modfc{k/2}{N}{\chi}$ 
the space of holomorphic modular forms of weight $k/2$, level $N$ and character
$\chi$.
We denote by
$\cusfc{k/2}{N}{\chi}$ 
and
$\eisfc{k/2}{N}{\chi}$ 
the subspace of cuspidal forms and the Eisenstein subspace, respectively.
When $\chi$ is the trivial character, we omit it from the notation.

We consider the following operators acting on half-integral weight modular
forms.
Let $g = \sum_{n \geq 0} a(n) q^n \in \modf{k/2}{N}$.
\begin{itemize}
	\item The Fricke involution $W(N)$, given by
\begin{align*}
	W(N) & : \modfc{k/2}{N}{\chi} \to \modfc{k/2}{N}{\chi \chi_N}, \\
		 & (g\vert W(N))(z) = (Nz)^{-k/2} g(-1/Nz).
\end{align*}
	We include here an extra factor of $N^{-k/2}$ not present in \cite{wang-pei}.
	\item For a prime $\ell$, the Hecke operator $T(\ell^2)$, given by
\begin{align}
	\nonumber
	T(\ell^2) & : \modfc{k/2}{N}{\chi} \to \modfc{k/2}{N}{\chi}, \\
	\label{eqn:hecke}
	g \Tls & = \sum_{n \geq 0} 
			  \left(a(\ell^2 n) 
				  + \chi(\ell) \ell^{\lambda-1} \omega_n(\ell) a(n)
		+ \chi(\ell^2) \, \ell^{2\lambda-1}a(n/\ell^2)\right)  q^n.
\end{align}
	\item For an integer $m \geq 1$, the $V(m)$ operator, given by
\begin{align*}
	V(m) & : \modfc{k/2}{N}{\chi} \to \modfc{k/2}{mN}{\chi \chi_m}, \\
		 & g \vert V(m) = \sum_{n \geq 0} a(n) q^{mn}.
\end{align*}
	\item For an integer $m \geq 1$, the $U(m)$ operator, given by
\begin{align*}
	U(m) & : \modfc{k/2}{N}{\chi} \to \modfc{k/2}{M}{\chi \chi_m},\\
		 & g \Um = \sum_{n \geq 0} a(mn) q^n,
\end{align*}
	where $M$ is the smallest multiple of $N$ which is divisible by every prime
	dividing $m$, and such that the conductor of $\chi_m$ divides $M$.
\end{itemize}
The latter two act as well on rings of formal power series.

\medskip

The following is the Sturm bound for general weights.
Its proof follows from the integral weight case; see \cite[Prop. 4.1]{rsst}.

\begin{prop}\label{prop:sturm}
	Let $k \geq 3$ be an integer, and let $m$ be a prime.
	Suppose that 
	$g=\sum_{n \geq 0} a(n) q^n \in \modf{k/2}{N} \cap \Zq$.
	Let
	\begin{equation*}
		n_0 = \left\lfloor\frac{k}{24} 
		\cdot [\SL_2(\Z):\Gamma_0(N)]\right\rfloor.
	\end{equation*}
	If $a(n) \equiv 0 \pmod{m}$ for $1 \leq n \leq n_0$, then 
	$g \equiv 0 \pmod{m\Zq}$.
\end{prop}

The result is also valid for proving equalities, namely when $m=0$.

\subsection*{Eta-quotients}

Let $\eta(z)$ denote the Dedekind eta function, which is given by
\begin{equation*}
	\eta(z) = q^{\tfrac1{24}}\prod_{n=1}^\infty \left(1-q^n\right),
	\qquad q=e^{2\pi iz}.
\end{equation*}
Given a finite set $X=\{(\delta,r_{\delta})\} \subseteq \Z_{>0}\times\Z$,
denote $s_X=\sum\delta r_{\delta}$.
Assuming that $s_X \equiv 0 \pmod{24}$, the eta-quotient defined by $X$ is
\begin{equation}
	\label{eqn:etaq}
	\eta^X(z)=\prod_{X}\eta(\delta z)^{r_{\delta}}
	= q^{\tfrac{s_X}{24}} \prod_X
	\prod_{n=1}^\infty \left(1-q^{\delta n}\right)^{r_\delta}
	\quad \in q^{\tfrac{s_X}{24}}\left(1+q\Zq\right).
\end{equation}
Note that $1/\eta^X$ is also an eta-quotient.

Let $k = \sum_X r_\delta$, and let $N$ be the smallest multiple of every
$\delta$, and of $4$ if $k$ is odd, such that
\[
	N \sum_X \frac{r_\delta}{\delta} \equiv 0 \pmod{24}.
\]
Finally, letting $m'=\prod_X \delta^{r_\delta}$ we let $m = m'$ for even
$k$, and $m = 2m'$ for odd $k$.
Then (see \cite[Thm. 3]{gordon-hughes} and \cite[Coro. 2.7]{treneer3}) we
have the following result.

\begin{prop}
	\label{prop:eta_wh}
	With the notation as above,
	$\eta^X$ is a weakly holomorphic modular form of weight $k/2$, level $N$ and
	character $\chi_m$.
\end{prop}

Thus, $\eta^X$ is holomorphic and nonzero in the upper half-plane, but it can
have poles and zeros at the cusps.
Furthermore, following \cite{ligozat}, if $\gcd(a,c) = 1$, then the
order of vanishing of $\eta^X$ at a cusp $s = a/c \in \Q\cup\{\infty\}$ is given
by
\begin{equation}
	\label{eqn:ligozat}
	\ord_s\left(\eta^X\right) =
	\frac{N}{24 \gcd(c^2,N)} \, \sum_X \gcd(c,\delta)^2 \, \frac{r_\delta}{\delta}.
\end{equation}

\begin{prop}
	\label{prop:Delta2}
	Let $\Delta_2 = \eta^8(z) \eta^8(2z)$.
	Then $\Delta_2 \in \cusf{8}{2}$.
	Furthermore, for every nonnegative integer $k$
	the map
	\begin{align*}
		\modf{k}{2} & \to \cusf{k+8}{2},\\
		g & \mapsto g\cdot \Delta_2
	\end{align*}
	is an isomorphism.
\end{prop}

\begin{proof}
	The above proposition gives that $\Delta_2 \in \modf82$.
	Furthermore, by \eqref{eqn:ligozat} we see that $\Delta_2$ has simple
	zeros at the cusps for $\Gamma_0(2)$, namely $0$ and $\infty$.
	Hence the second claim follows, since $\Delta_2$ does not vanish on the
	upper half-plane.
\end{proof}

\subsection*{Eisenstein spaces of integral weight and level 2}
We consider the subgroup 
$\Gamma_\infty = \set{\pm\smat{1 & n \\ 0 & 1}}{n \in \Z} \leq \SL_2(\Z)$.
Let $k \geq 2$ be an even integer.
Denote
\begin{align*}
	E_{k}(z) & = \sum_{\gamma \in \Gamma_\infty \backslash \SL_2(\Z)}
	\frac{1}{(c_\gamma z + d_\gamma)^k} \quad \in \modf{k}{1},\\
	D_2 & = 2 E_2 | V(2) - E_2 \quad \in \modf{2}{2}.
\end{align*}
Then $E_{k} \in 1 + q\Zq$.
Furthermore,
\begin{equation}
	\label{eqn:eisenstein}
	E_k = 1 - \frac{2k}{B_k} \sum_{n\geq1} \sigma_{k-1}(n) q^n.
\end{equation}

The following result will not be used in our proofs, but explains the type of
forms $h_m$ appearing in Table~\ref{tab:hm} below
(see also Remark~\ref{rmk:conjecture}).
Though it is probably well known, we give a proof for the sake of completeness.

\begin{prop}
	\label{prop:Mk2}
	Let $D_2,E_4$ be as above.
	Then for every nonnegative integer $k$ the set
	$\set{D_2^a E_4^b}{2a+4b = k}$ is a basis for $\modf{k}{2}$.
\end{prop}

\begin{proof}

	Denote by $\calV_k$ the subspace of $\modf{k}{2}$ generated by 
	$\set{D_2^a E_4^b}{2a+4b = k}$.
	Let $\Delta_2 = \eta^8(z) \eta^8(2z)$.
	Using Proposition~\ref{prop:sturm} and \eqref{eqn:eisenstein}
	we get that
	$$576\Delta_2 = 5 D_2^2 E_4 - E_4^2 - 4 D_2^4.$$
	Hence $\Delta_2 \in \calV_8$.
	Thus, to prove that $\modf{k}{2} = \calV_k$ for every $k$, by
	Proposition~\ref{prop:Delta2} it suffices to
	show that for every $f \in \modf{k+8}{2}$ there exists $g \in \calV_{k+8}$
	such that $f-g \in \cusf{k+8}{2}$.

	For this purpose it suffices to prove that there exist
	$g_\infty,g_0 \in \calV_{k+8}$ such that $g_\infty$ does not vanish at
	$\infty$ and $g_0$ vanishes at $\infty$ but not at $0$; equivalently,
	$g_0$ vanishes at $\infty$ but is not cuspidal.

	We can clearly let $g_\infty = D_2^a$ with $a = \tfrac{k+8}2$.
	In the case of $g_0$, it suffices to consider $k \in \{0,2,4,6\}$.
	Then using explicit bases for $\cusf{k+8}2$ we see that we can let $g_0$ be
	as in Table~\ref{tab:g0}.

	\begin{table}[ht]
	\begin{tabular}{c l}
		$k$ & $g_0$
		\\\hline
		$0$ & $D_2^4 - E_4^2$ \\
		$2$ & $D_2^5 - D_2 E_4^2$ \\
		$4$ & $D_2^6 - E_4^3$ \\
		$6$ & $D_2^7 - D_2 E_4^3$
	\end{tabular}
	\caption{Forms in $\calV_{k+8}$ vanishing at $\infty$ but not at $0$.  Used in the proof of Proposition~\ref{prop:Mk2}.}
	\label{tab:g0}
	\end{table}

	Finally, the independence of the forms $D_2^a E_4^b$ follows using the
	formulas for $\dim(\modf{k}{2})$ (see \cite{cohen-oesterle}).
\end{proof}

\section{Eisenstein spaces of half-integral weight and level 16}
\label{sect:eisen}

Wang and Pei (\cite{wang-pei}) considered the Eisenstein spaces of half-integral
weights, giving bases of eigenforms for these spaces in the case of level $4D$,
with $D$ odd and squarefree.
Relying on their definitions and results, we consider the case of level
$16$.
The main result of this section is the following.

\begin{prop}
	\label{prop:eis_eigen}
	Let $\ell \geq 3$ be prime,
	and let $k \geq 3$ be an odd integer.
	Then $T(\ell^2)$ acts by multiplication by $\sigma_{k-2}(\ell)$
	on $\eisf{k/2}{16}$.
\end{prop}

We also give in Proposition~\ref{prop:eis_coeffs} exact formulas for the
coefficients of the Eisenstein series, which are needed to prove the congruence
in \eqref{eqn:g11}.

\medskip

As in Section~\ref{sect:prelims},
let $\Gamma_\infty = \set{\pm\smat{1 & n \\ 0 & 1}}{n \in \Z}
\leq \SL_2(\Z)$.
Let $k \geq 3$ be an odd integer.
Denote $\lambda = \tfrac{k-1}2$.
Let $N \in \{4,8\}$.
For $\gamma \in \Gamma_0(N)$, let $j(\gamma,z)$ be the automorphy factor of
weight $1/2$.
For $k>3$ we denote
\begin{align*}
	E_{k,N}(z) & = \sum_{\gamma \in \Gamma_\infty \backslash \Gamma_0(N)}
	\frac{1}{j(\gamma,z)^k}, \\
	E'_{k,N} & = \tfrac{2^k N^\lambda}{1-(-1)^{\lambda}i} \cdot 
		E_{k,N} \vert W(N).
\end{align*}
For $k = 3$ we consider the difference
$E_{3,N} - 2 \sqrt N \, E'_{3,N}$ defined by the formulas above, which, for
simplicity, we will denote by $E_{3,N}$.

We start by giving the Fourier expansions of these Eisenstein series,
following \cite{wang-pei}.
For this purpose we introduce the following notation, which will not be used in
other parts of the article.

\medskip

For an even integer $v$ denote 
\begin{equation*}
	c_k^\pm(v) = \frac{1-2^{(2-k)v/2}}{1-2^{2-k}} \pm 2^{(2-k)v/2}.
\end{equation*}
Given a positive integer $n$, let $v_n = \val_2(n)$ and $n' = (-1)^\lambda n /
2^{v_n}$, and denote
\begin{align*}
	C_k(n) & =
	\begin{cases}
		c_k^-(v_n-1), & 2 \nmid v_n,\\
		c_k^-(v_n)  , & 2 \mid v_n,\, n' \equiv 3 \pmod 4, \\
		c_k^+(v_n) + 2^{((2-k)v_n + (3-k))/2} \, \kro{n'}2,
					  & 2 \mid v_n,\, n' \equiv 1 \pmod 4,
	\end{cases}
	\\
	\gamma_{k,4}(n) & =
	\begin{cases}
		C_k(n), & k > 3, \\
		C_3(n) - 2 , & k = 3,
	\end{cases}
	\\
	\gamma_{k,8}(n) & =
	\begin{cases}
		0 , & (-1)^\lambda n \equiv 2,3 \pmod 4, \\
		C_k(n) - 1, & (-1)^\lambda n \equiv 0,1 \pmod 4,\, k > 3, \\ 
		C_3(n) - 2, & (-1)^\lambda n \equiv 0,1 \pmod 4,\, k = 3.
	\end{cases}
\end{align*}

Let $\omega$ denote a Dirichlet character of conductor $f$.
Let $B_\lambda$ denote the $\lambda$-th Bernoulli polynomial.
Then we consider the generalized Bernoulli number
\begin{equation*}
	B_{\lambda,\omega} = 
	f^{\lambda - 1} \sum_{a = 1}^{f} \omega(a) B_\lambda(\tfrac a{f}).
\end{equation*}
Furthermore, letting $\mu$ denote the Möbius function, for each positive
integer $n$ we denote
\begin{align*}
	\beta_{\lambda,\omega}(n) & =
	\sum_{a,b} \mu(a)
	\omega(a)
	a^{-\lambda} b^{-2\lambda+1}, 
\end{align*}
where $a,b$ run over all positive odd integers such $(ab)^2 \mid n$.

Recall that for each positive integer $m$ we consider the primitive Dirichlet
character $\omega_m$ determined by 
\begin{align*}
	\omega_m(a) = \kro{(-1)^\lambda m}{a}, \quad \gcd(a,4m) = 1.
\end{align*}
We denote by $f_m$ its conductor, and we remark that ${f_m/m}$ is the square of
a rational number.
We let
\begin{equation}
	\label{eqn:alpha}
	\alpha_{\lambda,m} = 
	\frac
	{
		\sqrt{f_m/m} \,
		B_{\lambda,\omega_m}
	}
	{
		f_m^\lambda \,
		B_{2\lambda}}
	\frac{
		1-\omega_m(2) 2^{-\lambda}
	}
	{
		1-2^{-2\lambda}
	}.
\end{equation}

Finally, for each positive integer $n$ we denote
\begin{align}
	\label{eqn:an}
	a_{k,N}(n) & = 
			\alpha_{\lambda,n} \,
			\beta_{\lambda,\omega_n}(n) \,
			\gamma_{k,N}(n)
			\,
			n^\lambda,
			\\
	\label{eqn:apn}
	a'_{k,N}(n) & = 
			\alpha_{\lambda,nN} \,
			\beta_{\lambda,\omega_{nN}}(n)
			\,
			n^\lambda.
\end{align}

\begin{prop}
	\label{prop:eis_coeffs}
	For $N \in \{4,8\}$ and odd $k \geq 3$ we have
	\begin{align*}
		E_{k,N} & = 1 +
			\sum_{n \geq 1}
			a_{k,N}(n) \, q^n, \quad k \geq 3,\\
		E'_{k,N} & =
			\sum_{n \geq 1}
			a'_{k,N}(n) \, q^n, \quad k > 3.
	\end{align*}
\end{prop}

The proof is essentially given in \cite{wang-pei}; their formulas for the
coefficients of these Eisenstein series involve values $L(\lambda, \omega_m)$ of $L$-series of
quadratic characters at positive integers.
The latter are well known; we use them in the result below.

Given a positive integer $\lambda$ we denote
\begin{equation*}
	\label{eqn:elambda}
	e_\lambda = 
	\frac{
		2^{2+\lambda -k} \kro{2\lambda+3}2 \lambda!
	}
	{
		(1-2^{-2\lambda}) B_{2\lambda} \pi^\lambda
	}.
\end{equation*}

\begin{lemma}
	\label{lem:alphael}
	For every positive integer $m$ we have that
	\begin{equation*}
		\alpha_{\lambda,m} = e_\lambda \,
		\left(1-\omega_m(2) 2^{-\lambda}\right) \,
		L(\lambda,\omega_m) \,
		m^{-1/2}.
	\end{equation*}
	Moreover, $\sgn(\alpha_{\lambda,m}) = \kro{2\lambda+1}2$.
\end{lemma}

\begin{proof}
	From \cite[p. 337]{montgomery-vaughan} and \cite[Thm.
	9.17]{montgomery-vaughan} we have that for every quadratic character
	$\omega$ with conductor $f$ and such that $\omega(1) =
	(-1)^\lambda$ we have that
	\begin{equation*}
		L(\lambda,\omega) =
		\frac{
			\kro{2\lambda+3}2 2^{\lambda -1} \pi^\lambda \sqrt f
		}
		{
		{\lambda!}f^\lambda
		} \,
		B_{\lambda,\omega},
	\end{equation*}
	from which the first claim follows.

	The second claim follows from the fact that for every such $\omega$ we have
	that $L(\lambda,\omega) > 0$; hence
	\begin{equation*}
		\sgn(\alpha_{\lambda,m}) =
		\sgn(e_\lambda) =
		\kro{2\lambda+3}2 \sgn(B_{2\lambda}) =
		\kro{2\lambda+1}2.
		\qedhere
	\end{equation*}
\end{proof}

\begin{coro}
	\label{coro:notzero}
	Let $n$ be a squarefree positive integer.
	Then $a'_{k,N}(n) \neq 0$. Furthermore, $a_{k,N}(n) = 0$ if
	and only if $\gamma_{k,N}(n) = 0$.
\end{coro}

\begin{proof}[Proof of Proposition~\ref{prop:eis_coeffs}]
	Using the well known formulas for $\zeta(2\lambda)$ and
	$\Gamma(\lambda + 1/2)$, and using that
	\begin{equation*}
		\frac{
			(-i)^{\lambda+1/2} \left(1+(-1)^\lambda i\right)
		}
		{
			\sqrt 2
		}
		= \kro{2\lambda+1}2,
	\end{equation*}
	we obtain that
	\begin{equation*}
		e_\lambda = 
		\frac{
			(-2\pi i)^{\lambda+1/2} \left(1+(-1)^\lambda i\right)
		}
		{
			2^{2\lambda+1} \Gamma(\lambda+1/2)
			\zeta(2\lambda) (1-2^{-2\lambda})
		}
		.
	\end{equation*}
	Then the result follows straightforwardly from Lemma~\ref{lem:alphael} and
	the formulas \cite[(2.30), (2.33), (2.35), (2.36) and (2.38)] {wang-pei}.
\end{proof}

Proposition~\ref{prop:eis_coeffs} shows that $E_{k,N}$ and $E'_{k,N}$, which a
priori have their coefficients in a cyclotomic field (\cite[Thm. 2.3]{shimura}),
actually have rational coefficients.
The following results shows that, as in the integral weight case (see
\eqref{eqn:eisenstein}), their deno\-minators are controlled by $k$ and can be
described in terms of Bernoulli numbers.

\medskip

We will require the following result 
from Carlitz (\cite[Thms. 1 and 3]{carlitz}).

\begin{lemma}
	\label{lem:carlitz}
	Let $d$ be a fundamental discriminant,
	and let $\lambda$ be a positive integer.
	\begin{enumerate}
		\item If $d = -4$ and $\lambda$ is odd, then 
			$2 B_{\lambda,\chi_d}/\lambda \in \Z$.
		\item If $d = \pm p$, with $p>2$ prime, then
			$ B_{\lambda,\chi_d}/\lambda \in \Z_{(p)}$.
			Moreover, if 
			$2\lambda/(p-1)$ is an odd integer, then
			$p B_{\lambda,\chi_d} \in \Z$.
		\item Otherwise,
			$B_{\lambda,\chi_d}/\lambda \in \Z$.
	\end{enumerate}
\end{lemma}

We denote
\begin{equation*}
	S_\lambda = \begin{cases}
		2^{\val_p(\lambda)+1}, & 2\mid\lambda, \\
		1, & 2 \nmid \lambda.
		\end{cases}
\end{equation*}
Furthermore, we denote 
$S'_{\lambda,N} = S_\lambda$ (see Remark~\ref{rmk:betterSl} below).

\begin{prop}
	\label{prop:eis_denoms}
	For $N \in \{4,8\}$ and odd $k \geq 3$ we have
	\begin{align*}
		E_{k,N} & \in 1 + \tfrac{
			\lambda
			}{
			2^{\lambda-1}
			\left(2^{2\lambda}-1\right)
			B_{2\lambda}
			S_\lambda
		} 
			\, \Zq, \\
			E'_{k,N} & \in \tfrac{
				\lambda
				2^\lambda
			}{
				\left(2^{2\lambda}-1\right)
				B_{2\lambda}
				S'_{\lambda,N}
				N^\lambda
			}
			\,\Zq.
	\end{align*}
\end{prop}

\begin{proof}
	We prove the claim for $E_{k,N}$; the proof for $E'_{k,N}$ follows
	by similar arguments, using \eqref{eqn:apn}.

	Let $n$ be a positive integer.
	Recalling that $f_n$ denotes the conductor of $\omega_n$,
	write $n = f_n q_n^2 = f'_n (q'_n)^2$ with $f'_n$ squarefree, so that 
	$\sqrt{f_n / n} = 1 / q_n$ and
	$2q_n / q'_n \in \{1,2\}$.
	Moreover, let $w_n = \val_2(q'_n)$ and write
	$q'_n = 2^{w_n} q''_n$.
	Then letting
	\begin{align*}
		r_n & =
		{q''_n}^{2\lambda-1}\,
		\beta_{\lambda,\omega_n}(n)
		,
		\\
		s_n & =
		S_\lambda
		\left(2^\lambda-\omega_n(2)\right)
			B_{\lambda,\omega_n} / \lambda
		,
		\\
		t_n & =
		\left(2q_n/q''_n\right)^{2\lambda-1}\,
		\gamma_{k,N}(n)
		,
	\end{align*}
	and using \eqref{eqn:alpha}, according to \eqref{eqn:an} we can decompose
	\begin{equation*}
		a_{k,N}(n) =
		\frac{
			\lambda \,
		}
		{
			2^{\lambda-1} \,
			\left(2^{2\lambda}-1\right) 
			B_{2\lambda} \,
			S_\lambda
		}
		\cdot
			r_n \,
			s_n \,
			t_n
			\,.
	\end{equation*}

	From the definition of $\beta_{\lambda,\omega}(n)$ it is easy
	to see that $r_n \in \Z$.
	Furthermore, by the definition of $\gamma_{k,N}(n)$,
	we have that $t_n \in \Z$.
	To prove the result it suffices then to show that $s_n \in \Z$.

	First assume that $\lambda$ is odd and $n$ is a square.
	Then $s_n / S_\lambda = 2 B_{\lambda,\omega_n} / \lambda$,
	hence the claim follows by part (a) of Lemma~\ref{lem:carlitz}.

	Assume now that $\lambda$ is odd or $n$ is not a square.
	In case (c) of Lemma~\ref{lem:carlitz}, the claim follows immediately.
	In case (b), let $p = f_n$.
	Then the claim follows from quadratic reciprocity and
	\begin{equation}
		\label{eqn:congruence}
		\kro2p^{
		\frac{2\lambda}{p-1}
		}
		\equiv
		2^\lambda
		\pmod{p^{\val_p(\lambda)+1}},
	\end{equation}
	which holds for even $2\lambda/(p-1)$ as well.

	Finally, assume that $\lambda$ is even and $n$ is a square.
	Then $s_n = (2^\lambda-1) B_\lambda / \lambda$
	(unless $n=1$, when they differ by a sign).
	In this case the result follows from \eqref{eqn:congruence} and a result of
	Von Staudt, which asserts that the denominator of $B_\lambda / \lambda$
	equals
	\begin{equation*}
		\prod_{p-1 \mid \lambda} p^{\val_p(\lambda)+1}.
		\qedhere
	\end{equation*}
\end{proof}

\begin{rmk}
	\label{rmk:betterSl}
	Making considerations about the $2$-adic valuation of the generalized
	Bernoulli numbers, the result also holds letting
\begin{equation*}
	S_\lambda = 
	\begin{cases}
		1/2^{\lambda - 2}, & \text {for even } \lambda,\\
		1/2^{\lambda - 1}, & \text{for odd }\lambda,
	\end{cases}
	\qquad
	S'_{\lambda,4} = 
	\begin{cases}
		1/2, & \text{for }\lambda = 2,\\
		1/2^{\lambda + 1}, & \text{for even }\lambda > 2,\\
		1/2^{\lambda - 1}, & \text {for odd } \lambda,
	\end{cases}
\end{equation*}
and $S'_{\lambda,8} = 1/2^\lambda$.
Furthermore, the normalized Eisenstein series according to these sharper
constants seem to be primitive.

\end{rmk}

\begin{prop}
	\label{prop:eis_bases}
	Let $k \geq 3$ be odd.
	Then
	\begin{equation*}
		\dim \eisf{k/2}{16} =
		\begin{cases}
			4, & k = 3,\\
			6, & k > 3.
		\end{cases}
	\end{equation*}
	Furthermore,
	\begin{multline}
		\label{eqn:eis_gens}
		\eisf{k/2}{16} =  \\
		\begin{cases}
		\ll{
			E_{3,4},
			E_{3,4} \Vm[4],
			E_{3,8},
			E_{3,4}\Um[2] \Vm[2]
		}, & k = 3,\\
		\ll{
			E_{k,4},
			E_{k,4}\Vm[4],
			E'_{k,4},
			E'_{k,4}\Vm[4],
			E_{k,8},
			E'_{k,8}\Vm[2]
		}, & k > 3.
		\end{cases}
	\end{multline}
\end{prop}

\begin{proof}
	The first claim follows from \cite{cohen-oesterle}.

	Let $N \in \{4,8\}$.
	In \cite[Thm. 7.6]{wang-pei} it is proved that $E_{k,N} \in \eisf{k/2}{N}$.
	Considering the codomains of the operators $W(N),V(2),V(4),U(2)$ (see
	Section~\ref{sect:prelims}) we get that $\eisf{k/2}{16}$ contains the
	subspace on the right hand side of \eqref{eqn:eis_gens}, for $k \geq 3$.

	We now prove that the generators on the right hand side
	of \eqref{eqn:eis_gens} are linearly independent,
	using the formulas for their coefficients given by
	Proposition~\ref{prop:eis_coeffs}.
	Assume first that $k \equiv 5 \pmod 4$.
	Then
	\begin{align*}
		E_{k,4} & = 1 + 
			a_{k,4}(1) q + 
			a_{k,4}(2) q^2+ 
			a_{k,4}(3) q^3+ 
			a_{k,4}(4) q^4+ 
			a_{k,4}(5) q^5+ 
			O(q^6), \\
		E_{k,4} \Vm[4] & = 1 + 
			a_{k,4}(1) q^4 +
			O(q^6), \\
		E'_{k,4} & =
			a'_{k,4}(1) q + 
			a'_{k,4}(2) q^2+ 
			a'_{k,4}(3) q^3+ 
			a'_{k,4}(4) q^4+ 
			a'_{k,4}(5) q^5+ 
			O(q^6), \\
		E'_{k,4} \Vm[4] & =
			a'_{k,4}(1) q^4 +
			O(q^6), \\
		E_{k,8} & = 1 +
			a_{k,8}(1) q + 
			a_{k,8}(4) q^4+ 
			a_{k,8}(5) q^5+ 
			O(q^6), \\
		E'_{k,8} \Vm[2] & = 
			a'_{k,8}(1) q^2 + 
			a'_{k,8}(2) q^4+ 
			O(q^6).
	\end{align*}
	Then, since $a'_{k,4}(1) \, a'_{k,8}(1) \neq 0$ 
	(see Corollary~\ref{coro:notzero}),
	it suffices to prove that
	\begin{align*}
		\pmat{
			a_{k,4}(1) & a_{k,4}(3) & a_{k,4}(5) \\
			a'_{k,4}(1) & a'_{k,4}(3) & a'_{k,4}(5) \\
			a_{k,8}(1) & 0 & a_{k,8}(5)
		}
	\end{align*}
	is non-singular.

	We have that $\beta_{\lambda,\omega}(n) = 1$ for squarefree $n$.
	Furthermore, we have that $\gamma_{k,4}(1) > 0, \gamma_{k,4}(3) < 0,
	\gamma_{k,4}(5) >0$ and that $\gamma_{k,8}(1) > 0, \gamma_{k,8}(5) < 0$.
	Then by Lemma~\ref{lem:alphael} the signs of the matrix above are given by
	\begin{align*}
		\kro{2\lambda+1}2 
		\pmat{
		 	+ & - & + \\
			+ & + & + \\
			+ & 0 & -
		},
	\end{align*}
	hence its determinant is non-zero.
	
	The case $k \equiv 7 \pmod 4, k > 3$, can be proved similarly, using the
	$7$-th coefficient instead of the $5$-th coefficient in the matrix above.
	Finally, for $k = 3$ using Proposition~\ref{prop:eis_coeffs} we get that
	\begin{align*}
		E_{3,4} & = 1 + 6q + 12q^2 + 8 q^3 + O(q^4), \\
		E_{3,4} \Vm[4] & = 1 + O(q^4), \\
		E_{3,8} & = 1 + 8q^3 + O(q^4), \\
		E_{3,4} \Vm[2] \Um[2] & = 1 + 12q^2 + O(q^4),
	\end{align*}
	which completes the proof.
\end{proof}

\begin{proof}[Proof of Proposition~\ref{prop:eis_eigen}]
	Denote by $\calV \subseteq \eisf{k/2}{16}$ the
	$\sigma_{k-2}(\ell)$-eigenspace for $T(\ell^2)$.

	We claim first that $E_{k,4}, E_{k,8} \in \calV$.
	For every $n$ we see easily from the definitions and Lemma~\ref{lem:alphael}
	that
	\begin{align*}
		\omega_{\ell^2 n} & = \omega_{n}, \\
		\alpha_{\lambda,\ell^2 n} & = \ell^{-1} \, \alpha_{\lambda,n}, \\
		\gamma_{k,N}\left(\ell^2 n\right) & = \gamma_{k,N}(n).
	\end{align*}
	Then the claim follows directly from \eqref{eqn:hecke},
	using the equalities above and the transformation
	formulas for computing
	$\beta_{\lambda,\omega_{\ell^2 n}}(\ell^2 n)$ in terms of
	$\beta_{\lambda,\omega_n}(n)$ given in \cite[p. 209]{wang-pei};
	we remark that though Wang and Pei are considering $k > 3$ and level $4D$
	with $D$ odd and squarefree, these particular computations hold in our
	setting.

	The result follows then by noting that the remaining generators for
	$\eisf{k/2}{16}$ given in Proposition~\ref{prop:eis_bases} belong to
	$\calV$, since by \cite[Thm. 5.19]{wang-pei} the
	Hecke operators $T(\ell^2)$ with $\ell \neq 2$
	commute with the operators $W(N)$,
	and by $\eqref{eqn:hecke}$ they commute with $U(2),V(2),V(4)$.
\end{proof}

\section{Proofs}
\label{sect:proofs}

This section is devoted to give the proofs of the theorems stated in the
Introduction.

\medskip

We first state the following result for obtaining congruences for
coefficients of (modulo $m$) eigenforms of half-integral weight, used by
\cite{ahlgrenono,treneer2,treneer1,rsst}  among others.

\begin{prop}
	\label{prop:congruences}
	Let $g = \sum_{n\geq 0} a(n) q^ n \in \modf{k/2}{N}
	\cap \Z\llbracket q\rrbracket$,
	and let $\ell,m$ be primes such that 
	$g \vert T(\ell^2) \equiv \eigen \, g \pmod{m\Zq}$.
	\begin{enumerate}
		\item If $\eigen \equiv 0 \pmod m$, then
			$a(\ell^3 n) \equiv 0 \pmod m$ 
			for every $n$ prime to $\ell$.
		\item If there exists $\epsilon \in \{\pm1\}$ such that 
			\begin{equation*}
				\eigen
				\equiv 
				\epsilon \, \ell^{\tfrac{k-3}2} \pmod m,
			\end{equation*}
			then $a(\ell^2 n) \equiv 0 \pmod m$ 
			for every $n$ prime to $\ell$ such that
			$\omega_n(\ell) = \epsilon$.
	\end{enumerate}
\end{prop}

\begin{proof}
	Both claims follow directly from \eqref{eqn:hecke}; for part (a), by
	replacing $n$ by $\ell n$, with $n$ prime to $\ell$.
\end{proof}

The goal of the following series of results is to prove that for every odd prime $m$ the numbers
$\overline{p}(mn)$ are congruent modulo $m$ to the Fourier coefficients of a holomorphic
modular form.
We start with two preliminary results.

\begin{lemma}
	\label{lem:Power}
	Let $f$ and $g$ be power series, and let $m \geq 1$.
	Then
	\begin{equation*}
		(\left(f\Vm \cdot g\right)\Um = f\cdot \left(g\Um\right).
	\end{equation*}
\end{lemma}

\begin{proof}
	Let $f=\sum_{n=0}^{\infty}a(n)q^n$ and
	$g=\sum_{n=0}^{\infty}b(n)q^n$. 
	Denote
	\begin{align*}
		\widetilde{a}(h) & =
			\begin{cases}
			a(n), &\text{if $h=nm$},\\
			0, &\text{otherwise},
			\end{cases}
			\\
			\widetilde{c}(h) & 
			=\sum_{k=0}^{h}\widetilde{a}(k)b(h-k).
	\end{align*}
	Now, note that
	\begin{align*}
		\widetilde{c}(hm)=\sum_{k=0}^{hm}\widetilde{a}(k)b(hm-k)
			=\sum_{k=0}^{h}a(k)b(hm-km).
	\end{align*} 
Then we have
	\begin{align*}
		((f \Vm) \cdot g)\Um&=
		\left(
		\left(\sum_{h=0}^{\infty}\widetilde{a}(h)q^{h}\right)
		\left(\sum_{n=0}^{\infty}b(n)q^n\right)
		\right)\Um\\
				&=\left(\sum_{h=0}^{\infty}\widetilde{c}(h)q^{h}\right)\Um
				=\sum_{h=0}^{\infty}\widetilde{c}(hm)q^{h} \\
				 & =\sum_{h=0}^{\infty}
				\left(\sum_{k=0}^{h}a(k)b(mh-k)\right) q^h = f\cdot \left(g\Um\right).
		\qedhere
	\end{align*} 
\end{proof}

\begin{lemma}
	\label{lem:Modm}
	Let $f$ be an eta-quotient.
	Then for every prime $m \geq 1$ we have that
	\begin{equation*}
		f\Vm\equiv f^{m}\pmod{m\Zq}.
	\end{equation*}
\end{lemma}

\begin{proof}
	Write $f$ as in \eqref{eqn:etaq}.
	Since both operators $V(m)$ and $g\mapsto g^m$ are multiplicative, it
	suffices to verify the congruence for every factor $g$ of $f$.

	For $g = q^{\tfrac{s_X}{24}}$ both operators clearly agree, and for 
	$g = 1-q^{\delta n}$ the congruence follows from the fact that 
	$(r+s)^m \equiv r^m + s^m \pmod{m\Zq}$ for every $r,s \in \Zq$.
\end{proof}

In what follows we consider the eta-quotient related to the generating function for
$\overline{p}(n)$ (see \cite[(1.1)]{overpartitions}).
Namely, we let
\begin{equation}
	\label{eqn:etaq_op}
	f = \frac{\eta(2z)}{\eta^2(z)}
	= \sum_{n\geq 0} \overline{p}(n)q^n,
\end{equation}
We remark that $f$ is not holomorphic: by \eqref{eqn:ligozat}, it has a simple
pole at $s=0$.

\begin{lemma}
	\label{lem:F16}
	Denote $F = 1/f$.
	Then $F \in \modf{1/2}{16}$.
\end{lemma}

\begin{proof}
	By Proposition~\ref{prop:eta_wh} we have that $F$ is a weakly holomorphic
	modular form of level $16$ and weight $1/2$, with trivial character.
	Its possible singularities lie at the cusps $s$ for $\Gamma_0(16)$, namely $s
	\in \{0,1/8,1/4,1/2,3/4,\infty\}$.
	Then the claim follows from \eqref{eqn:ligozat}, which shows
	that the order of vanishing of $F$ at each such $s$ is nonnegative
	(moreover, it is positive only for $s=0$).
\end{proof}

For the following two propositions we let $0 < a_m < 8$ be such that 
$a_m \equiv -m \pmod 8$, and we denote
\begin{equation*}
	r_m = \tfrac12 \left(m(16-a_m) -17\right).
\end{equation*}

\begin{prop}
	\label{prop:hm'}
	Let $m \geq 3$ be a prime.
	There exists $h'_m \in \modf{r_m}{2} \cap \Zq$
	such that
	\begin{equation}
		\label{eqn:fUmgm'}
		f \Um \equiv F^{a_m} h'_m \pmod {m\Zq}.
	\end{equation}
	In particular, $f \Um$ is congruent modulo $m \Zq$ to a holomorphic modular
	form.
\end{prop}

\begin{proof}

	Recall the eta-quotient $\Delta_2(z) = \eta^8(z) \eta^8(2z)$.
	We consider the eta-quotients
	\begin{equation*}
		\alpha = \Delta_2 F^{-a_m}, \qquad
		\beta = f \alpha^m.
	\end{equation*}
	By \eqref{eqn:ligozat} we have that $\beta \in \cusf{r_m+8}{2}$.
	Since Hecke operators preserve cuspforms, by Proposition~\ref{prop:Delta2}
	there exists $h'_m \in \modf{r_m}{2}$ such that
	\begin{equation*}
		\beta \vert T(m) = \Delta_2 h'_m.
	\end{equation*}
	Note that since $\beta \in q\Zq$ and $\Delta_2 \in q\Zq^\times$ we have that $h'_m
	\in \Zq$.

	On the other hand, using Lemmas~\ref{lem:Power} and \ref{lem:Modm} we have that
	\begin{equation*}
		\beta \Um \equiv
		(f \cdot (\alpha \Vm)) \Um \equiv
		f \Um \cdot \alpha
		\pmod{m\Zq}.
	\end{equation*}
	Since over integral weights $T(m)$ agrees with $U(m)$ modulo $m\Zq$,
	the above congruences give that
	\begin{equation*}
		f \Um \cdot \alpha \equiv \Delta_2 h'_m
		\pmod{m\Zq},
	\end{equation*}
	which, since $\alpha,\Delta_2 \in q\Zq^\times$, concludes the proof.
\end{proof}

Proposition~\ref{prop:hm'}, together with the following refinement of
\cite[Thm 1.1]{treneer2} imply Theorem~\ref{thm:treneer} straightforwardly.

\begin{thm}[Treneer]

	Let $f = \sum_n a(n) q^n$ be a weakly holomorphic modular form with level $N$
	and integral coefficients.
	Let $m$ be an odd prime such that $f \Um$ is congruent modulo $m \Zq$ to a
	holomorphic modular form.
	Then a positive proportion of the primes $\ell \equiv -1 \pmod{N m}$ have
	the property that
	\begin{align}
		\label{eqn:treneer}
		a\left(m \ell^3 n\right) \equiv 0  \pmod m.
	\end{align}
	for every $n$ coprime to $m\ell$.
\end{thm}

\begin{proof}
	The proof given by Treneer, in which a weaker version of
	\eqref{eqn:treneer} is obtained, holds with no changes;
	by introducing the hypothesis on $f \Um$ we get that its main ingredient,
	namely \cite[Prop. 3.5]{treneer2}, holds for $m$ instead of for a large enough
	power of $m$.
\end{proof}

We now show that, at least for small values of $m$, the weight of the
holomorphic modular form in Proposition~\ref{prop:hm'} can be improved.

\begin{prop}
	\label{prop:hm}
	Let $3 \leq m \leq 19$ be a prime.
	Let $h_m$ be the corresponding form given in Table~\ref{tab:hm},
	and let $g_m = F^{a_m} h_m$.
	Then $g_m \in \modf{k_m/2}{16} \cap \Zq$, and
	\begin{equation}
		\label{eqn:fUmgm}
		f \Um \equiv g_m \pmod {m\Zq}.
	\end{equation}
\end{prop}

\begin{table}[ht]
\begin{tabular}{c l}
	$m$ & $h_m$
    \\\hline
	$3$ & $1$ \\
	$5$ & $1$ \\
	$7$ & $D_2$ \\
	$11$ & $D_2$ \\
	$13$ & $E_4$ \\
	$17$ & $13 D_{2}^{2} + 5 E_{4}$\\
	$19$ & $11 D_{2}^{3} + 9 D_{2} E_{4}$ \\
\end{tabular}
\caption{Holomorphic modular forms used in Proposition~\ref{prop:hm}.}
\label{tab:hm}
\end{table}

\begin{proof}
	The first claim follows from Lemma~\ref{lem:F16},
	since for every $m$ we have that $h_m \in \modf{\tfrac{k_m-a_m}{2}}{2}$.
	We have to verify \eqref{eqn:fUmgm}.

	Assume first that $m>3$.
	Let $e = (15-a_m)/2$.
	Since $E_{m-1} \equiv 1 \pmod{m\Zq}$, it suffices to prove that
	the form $h'_m$ from the above proposition satisfies that
	\begin{equation*}
		h'_m \equiv h_m E_{m-1}^e
		\pmod{m\Zq}.
	\end{equation*}
	With the above choice of $e$, both forms in this congruence belong to
	$\modf{r_m}{2} \cap \Zq$.
	Thus, by Proposition~\ref{prop:sturm} and \eqref{eqn:fUmgm'}
	it suffices to prove that the $n$-th
	coefficients of $f\Um \cdot f^{a_m}$ and $h_m$ agree, modulo $m$, up to $n$
	equal to $n_0 = \left\lfloor\tfrac{r_m}{36}\right\rfloor$.
	In each case, this can be proved by computing these numbers explicitly.

	The case $m=3$ follows similarly, replacing $E_{m-1}^e$ by $D_2^4$.
\end{proof}

\begin{rmk}
	\label{rmk:conjecture}
	In fact, using the techniques from the above proof and
	Proposition~\ref{prop:Mk2}, we have found forms
	$h_m$ as in Proposition~\ref{prop:hm} for every prime
	$m < 1000$.
\end{rmk}

\color{black}

From here on, given primes $m,\ell$, we denote
\[
	\lambda_{m,\ell} = 1 + \ell^{k_m-2},
\]
the eigenvalue of $T(\ell^2)$ acting on $\eisf{k_m/2}{16}$ 
(see Proposition~\ref{prop:eis_eigen}).

\begin{prop}
	\label{prop:gm_eigens}
	Let $3 \leq m \leq 19$ be a prime,
	and let $g_m$ be the form given in Proposition~\ref{prop:hm}.
	\begin{enumerate}
		\item If $3 \leq m \leq 11$ then $g_m \Tls \equiv \eigen g_m
			\pmod{m\Zq}$ for every prime $\ell > 2$.
		\item If $13 \leq m \leq 19$ then 
			$g_m \Tls \equiv \eigen g_m \pmod{m\Zq}$
			for every prime $\ell$ in Table~\ref{tab:eigens}.
	\end{enumerate}
\end{prop}

\begin{table}[ht]
\begin{tabular}{c|l}
	$m$ & $\ell$
    \\\hline
	$13$ & $431, 1811, 1871, 1949, 2207, 2459, 3301, 4001, 4079, 4289, 4513, 4799, 4931$ \\
	$17$ & $1999, 2207, 2243, 4759$\\
	$19$ & $151, 1091, 2207, 2659, 3989$ \\
\end{tabular}
\caption{Primes $\ell$ giving congruences modulo $m$.  See Proposition~\ref{prop:gm_eigens}.}
\label{tab:eigens}
\end{table}

\begin{proof}
	To prove part (a) we can use Proposition~\ref{prop:eis_eigen}, once we
	verify that for $3 \leq m \leq 11$ we have that 
	$g_m \in \eisf{k_m/2}{16} + m\Zq$.
	The latter claim, in the case $3 \leq m \leq 7$, follows from the fact that
	$\cusf{k_m/2}{16} = \{0\}$. 
	In the case $m = 11$, since 
	by Proposition~\ref{prop:eis_denoms} we have that
	\begin{equation*}
		2^3 17 \cdot E_{9,4},\; 
		2^3 17 \cdot E_{9,8},\; 
		2^4 17 \cdot E'_{9,4},\;
		2^5 17 \cdot E'_{9,8} 
		\quad \in \Zq,
	\end{equation*}
	we can use Proposition~\ref{prop:sturm} to obtain that
	\begin{multline}
		\label{eqn:g11}
		g_{11}(z) \equiv
					9 E_{9,4} + 
					4 E_{9,4} \vert V(4) +
					7 E'_{9,4} + 
					4 E'_{9,4} \vert V(4) +
					7 E'_{9,8} \vert V(2)
					\pmod{11\,\Zq}.
	\end{multline}

	In order to prove part (b), by Proposition~\ref{prop:sturm} it suffices to prove
	that the $n$-th coefficients of $g_m\Tls$ and $\eigen \, g_m$ agree,
	modulo $m$, for $n$ up to
	\[
		\tfrac{k_m}{24} \cdot \left[\SL_2(\Z): \Gamma_0(16)\right]
			= m-2.
	\]
	Moreover, by \eqref{eqn:fUmgm} it suffices to prove that
	\[
		\left(f\Um\Tls\right)(n) \equiv 
		\eigen \, (f\Um)(n) \pmod m,
		\quad 1 \leq n \leq m-2,
	\]
	which in each case can be proved by computing these numbers explicitly.
\end{proof}

\begin{rmk}
	The proof of Proposition~\ref{prop:hm} involves computing $\pb(mn)$ 
	modulo $m$ for small values of $n$.
	This can be accomplished easily by expanding the infinite product
	\eqref{eqn:etaq_op} defining $f$.

	The proof of Proposition~\ref{prop:gm_eigens} involves computing 
	$\pb(mn)$ modulo $m$ for large values of $n$ (e.g. $n =
	m(m-2)\ell^2$ with large $\ell$); in this case we resort to the
	efficient method provided by \cite{barquero-etal}.
\end{rmk}

The proofs of our main results now follow easily.

\begin{proof}[Proof of Theorems~\ref{thm:inf0} and~\ref{thm:fin0}]
	They follow using Pro\-position~\ref{prop:congruences} (a) and
	Proposition \ref{prop:gm_eigens}.
\end{proof}

\begin{proof}[Proof of Theorems~\ref{thm:infeps} and~\ref{thm:fineps}]
	They follow using Pro\-position~\ref{prop:congruences} (b) and
	Proposition \ref{prop:gm_eigens};
	the eigenvalues $\eigen$ in Table~\ref{tab:eigens} satisfy the
	hypothesis of Proposition~\ref{prop:congruences} (b), namely they are such
	that
	\begin{equation*}
		\eigen \equiv \eml \, \ell^{\tfrac{k_m-3}2}\pmod m,
	\end{equation*}
	where $\eml$ is as in Table~\ref{tab:ml_eps}.
\end{proof}

\begin{rmk}
	We found that $g_m$ is, modulo $m\Zq$, an eigenfunction of
	$T(\ell^2)$ for more primes $\ell$ than those appearing in
	Table~\ref{tab:eigens}, but the eigenvalues are not useful for our purposes,
	since they do not satisfy any of the hypotheses of
	Proposition~\ref{prop:congruences}.
	Moreover, the primes given are all the primes $\ell < 5000$ giving
	congruences.
	
	For $m = 23$ we found that $\ell = 5303, 8783$ yield eigenvalues, but
	they do not give congruences.
	For larger $m$ we have not been able to find eigenvalues.
\end{rmk}


\end{document}